 \documentclass[11pt]{article}
\usepackage{amssymb,amsmath,latexsym}
\usepackage{epsfig}
\usepackage{enumerate}

\allowdisplaybreaks

\newtheorem{thm}{Theorem}[section]
\newtheorem{cor}[thm]{Corollary}
\newtheorem{lemma}[thm]{Lemma}

\newtheorem{remark}[thm]{Remark}
\newtheorem{example}[thm]{Example}
\numberwithin{equation}{section}

\def\pf{{\medskip\noindent {\bf Proof. }}}
\def\qed{{\hfill $\Box$ \bigskip}}

\def\sD {{\cal D}}  
  
  \def\sL {{\cal L}}
 \def\sN {{\cal N}}

 \def\bB {{\mathbb B}} 
 \def\bE {{\mathbb E}}

  \def\bR {{\mathbb R}}

\def\E{{\mathbb E}}
\def\P{{\mathbb P}}

\def\bea{\begin{align*}}
\def\eea{\end{align*}}
\def\bee{\begin{equation}}
\def\eee{\end{equation}}

\def\eps{\varepsilon}

\makeatletter
\@addtoreset{equation}{section}

\makeatother

\begin{document}
\bibliographystyle{plain}

\title{\Large {\bf Time fractional  equations and probabilistic representation} \medskip \\
\small (In {\it  Chaos, Solitons and Fractals \bf 102} (2017), 168-174.)}

\author{{\bf Zhen-Qing Chen} }
\date{January 12, 2019}
\maketitle

\begin{abstract}
In this paper, we study the existence and uniqueness of solutions for general  fractional-time parabolic equations  of mixture type,
 and their  probabilistic representations   in terms of  the corresponding inverse subordinators with or without drifts. 
An explicit relation between occupation measure for Markov processes
time-changed by inverse subordinator in open sets and that of the original Markov process in the  open set is also given. 
\end{abstract}



\bigskip

\section{Introduction}

Fractional calculus has attracted lots of attentions in several fields including mathematics, physics, chemistry, engineering, hydrology
 and even finance and social sciences (see \cite{H11,MSi, MK,MK04}).
The classical heat equation $\partial_t u=\Delta u$ describes heat propagation in homogeneous medium.
The time-fractional diffusion equation $\partial^\beta_t u=\Delta  u$
with $0<\beta<1$  has been widely used
to model the anomalous diffusions  exhibiting  subdiffusive behavior,
due to particle sticking and trapping phenomena (see e.g. \cite{MSi,SKW82}).
Here  the fractional-time derivative $\partial^{\beta}_t$ is the Caputo derivative of order $\beta\in (0,1)$, which can be defined by
\begin{equation}\label{e:1.1}
\partial_t^\beta f(t)=\frac{1}{\Gamma (1-\beta)} \frac{d}{dt}
\int_0^t (t-s)^{-\beta}  \left(f(s)-f(0)\right) ds , 
\end{equation}
where $\Gamma(\lambda):=
\int^{\infty}_0 t^{\lambda-1} e^{-t}dt$  is the Gamma function. The above definition says that the fractional
derivative of $f$ at  time $t$ depends on the whole history of $f(s)$ on $(0, t)$ with  the nearest past affecting the present more. 
 Meerschaert and Scheffer \cite[Theorem 5.1]{MS0} recognized, 
  based on Baeumer and Meerschaert \cite{BMe}, that the solution to
$u=u(t, x)$ of  $\partial^\beta_t u=  \Delta  u$ with $u(0, x)=f(x)$
admits an interesting probabilistic representation: 
$$
u(t, x)=\E_x[ f(X_{E_t})], \quad  x\in \bR^d,
$$
where $X$ is Brownian motion on $\bR^d$ with infinitesimal generator $\Delta$ 
and $E_t$ is an inverse 
 $\beta$-stable subordinator that is independent of $X$. 
In fact, the above representation was proved for a large class of operators $\sL$ 
in place of $\Delta$   that generates a strong Markov process $X$. 
This representation  connects   probability theory to  time fractional equations.
 The scaling property of the $\beta$-stable subordinator is used in a crucial way in their derivation.
 
In applications and numerical approximations \cite{DYZ}, there is  a need to consider generalized fractional-time derivatives  where 
its value at time $t$   may depend  only on the finite range of the past from $t-\delta$ to $t$, for example, 
$  \frac{d}{dt}
\int_{(t-\delta)^+}^t (t-s)^{-\beta}  \left(f(s)-f(0)\right) ds$. Here for $a\in \bR$, $a^+:=\max\{a, 0\}$. 
Motivated by this, for a given   function $w: (0, \infty) \to [0, \infty)$ that is locally integrable on $[0, \infty)$, 
we   introduce a generalized fractional-time derivative
\begin{equation}\label{e:1.2}
\partial_t^w f(t)= \frac{d}{dt}
\int_0^t w(t-s)  \left(f(s)-f(0)\right) ds,
\end{equation} 
whenever it is well defined. 
Typically $w(t)$ is a non-negative decreasing function on $(0, \infty)$ that blows up at $t=0$.
Clearly, when $w(s)=\frac1{\Gamma (1-\beta)} s^{-\beta}$ for   $\beta \in (0, 1)$, $\partial^w_tf$ is just the Caputo derivative of order $\beta $
defined by \eqref{e:1.1}. 

  Let $X=\{X_t, t\geq 0; \, \P_x, x\in E\}$ be a strong Markov process on a separable locally compact Hausdorff space $E$ 
  whose transition semigroup $\{P_t, t\geq 0\}$ is a uniformly bounded strong continuous semigroup in some Banach space $(\bB, \| \cdot \|)$.
For example,     $\bB=L^p(E; m)$ for some measure $m$ on $E$ and $p\geq 1$ or $\bB=
   C_\infty (E)$, the space of continuous functions  on $E$ that vanish at infinity equipped with uniform norm.
  Let  $(\sL, \sD (\sL))$ be the infinitesimal generator of $\{ P_t, t\geq 0\}$ in $\bB$.  
	In this paper, we are interested in  the existence and uniqueness
  of solution $u=u(t, x)$  for
$$ 
\kappa  \frac{\partial u}{\partial t}+ \partial_t^w u= \sL u \quad \hbox{with } u(0, x)=f(x)  
$$
  and its probabilistic representation, where $\kappa \geq 0$ is a positive constant. 
We will also address the following question: given a subordinator $S=\{S_t; t\geq 0\}$ that is independent of $X$,
what equation does $u(t, x):=\E_x \left[ f(X_{E_t})\right]$ satisfy?

Given a constant $\kappa \geq 0$ and an unbounded right continuous non-increasing function $w(x)$ on $(0, \infty)$
with $\lim_{x\to \infty } w(x)=0$ and $\int_0^\infty (1\wedge x) (-dw(x))<\infty$, there is a unique 
non-negative valued L\'evy process $\{S_t; t\geq 0\}$ with $S_0=0$ (called subordinator) associated with it in the following way.
Here for $a, b\in \bR$, $a\wedge b:=\min\{a, b\}$. 
Let $\mu$ be the measure on $(0, \infty)$ so that $w(x)=\mu (x, \infty)$. Clearly 
$$
\mu (0, \infty)=\infty \quad \hbox{and} \quad  \int_0^\infty (1\wedge x) \mu (dx)<\infty.
$$
It is well-known (cf. \cite{Be}) that there is subordinator $\{S_t; t\geq 0\}$   with
    Laplace exponent $\phi$:
  \begin{equation}\label{e:1.3}
   \E \left[ e^{-\lambda S_t} \right] = e^{-t \phi (\lambda)}, \quad \lambda >0,
 \end{equation}
 so that  
 \begin{equation}\label{e:1.4} 
 \phi (\lambda) =\kappa \lambda + \int_0^\infty (1-e^{-\lambda x}) \mu ( dx).
 \end{equation}
 The measure $\mu$ is called the L\'evy measure of the subordinator.
 
 Conversely,  given a subordinator  $\{S_t; t\geq 0\}$, 
 there is a unique constant $\kappa\geq 0$ and a L\'evy measure $\mu$  on
$(0, \infty)$ satisfying 
 $\int_0^\infty (1\wedge x) \mu (dx)<\infty$ so that \eqref{e:1.3} and \eqref{e:1.4} hold. 
Throughout this paper, $\{S_t; t\geq 0\}$ is such a general subordinator   with infinite L\'evy measure $\mu$
and  possibly with drift $\kappa \geq 0$.  When $\kappa =0$, we say the subordinator is driftless or with no drift. 
   Define for $t>0$, 
$E_t =\inf\left\{ s>0: S_s >t\right\}$, the inverse subordinator. 
 The assumption that the L\'evy measure $\mu$ is infinite (which is equivalent to 
 $w(x):=\mu (x, \infty)$ being unbounded) excludes  compounded Poisson processes. 
  Under this assumption, almost surely, $t\mapsto S_t$ is strictly increasing and hence $t\mapsto E_t$ is continuous. 
 
The main purpose of this paper is to establish the following.

\begin{thm}\label{T:1.1}   Under the above setting, let $w(x)= \mu  (x, \infty)$, 
which is an unbounded right continuous non-increasing function on $(0, \infty)$.  
The function $u(t, x):=\E_x [ f(X_{E_t})]$ is the unique solution in $\bB$ to the time fractional equation 
 \begin{equation}\label{e:1.5}
\left( \kappa    \partial_t   + \partial_t^w\right)u= \sL u \quad \hbox{with } u(0, x)=f(x) 
 \end{equation}
 in the strong sense (see Theorem \ref{T:2.3} for a precise statemnt)  for every $f\in \sD (\sL)$ . Here $\partial_t$ is the time derivative $\frac{\partial}{\partial t}$.  
\end{thm}

Our method of  proof to the above theorem is  different from that of \cite{BMe} which is for stable subordinators,  
as there is no scaling property for a general subordinator $S_t$. 
Our approach is  quite robust and direct that works for any subordinator with infinite L\'evy measure  
 and for a wide class of infinitesimal generators.  One feature of this paper is that possible  mixture of the standard time derivative $\partial_t$
 and the general fractional time derivative $\partial^w_t$ is covered and treated in a unified way.  
Moreover, we will establish a more general result for $\sL$ being the infinitesimal generator of any uniformly bounded 
strongly continuous semigroup  in general Banach spaces;   see Theorem \ref{T:2.3} for a  precise statement. 
Our Theorem \ref{T:2.3} not only gives the existence but also the uniqueness of  solutions to the time fractional equation.  
The generalized Caputo derivative defined by \eqref{e:1.2} with  $w(x)=\mu (x, \infty)$ extends the 
distributed order fractional derivative defined in \cite{MS1} where $S_t$ is a   mixture of $\beta$-stable subordinators. 
An important application of these more general time fractional derivatives 
is to model ``ultraslow diffusion" where a plume spreads at a logarithmic rate; see \cite{MS1}
for details.
 
  In Section 3 of this paper, we will study the relation between occupation measure for
the time-changed process $X^*:=X_{E_t}$ by inverse subordinator in an open set $D\subset E$ with that of $X$ in $D$. 

\section{General time fractional equations}\label{S:2}

Recall that  $\{S_t; t\geq 0\}$ is a general subordinator with infinite L\'evy measure $\mu$ and drift $\kappa \geq 0$, whose Laplace exponent $\phi (\lambda)$ is given by \eqref{e:1.4}.
Define $w(x)=\mu (x, \infty)$ for $x>0$ and $\phi_0(\lambda ) :=\int_0^\infty \left( 1- e^{-\lambda x} \right) \mu (dx)$.
Note that $\phi_0 (\lambda)$ is the Laplace exponent of
 the driftless subordinator  $\{\bar  S_t:=S_t-\kappa t, t\geq 0\}$ having L\'evy measure $\mu$. 
Clearly
\begin{equation}\label{e:2.1}
\phi (\lambda)= \kappa \lambda + \phi_0 (\lambda) 
\quad \hbox{ and } \quad S_t=\kappa t + \bar S_t.
\end{equation}
Since $\mu (0, \infty)=\infty$, almost surely, $t\mapsto \bar S_t$ is strictly increasing. 

For every $a>0$, by Fubini theorem, 
\begin{equation}\label{e:2.2}
\int_{0 }^a   w(x) dx  =\int_{0 }^a  \left( \int_{(x, \infty)}  \mu (d\xi) \right)  dx
=\int_{0 }^\infty \left( \int_0^{\xi \wedge a} dx \right) \mu (d \xi)=
 \int_0^\infty
(\xi \wedge a ) \mu (d \xi)<\infty.
\end{equation}

The Laplace transform of $w$ is
\begin{eqnarray}
\int_0^\infty e^{-\lambda x} w(x) dx
&=& \int_0^\infty e^{-\lambda x} \int_{(x, \infty)} \mu (d\xi) dx   
 =  \int_{0 }^\infty \left( \int_0^\xi e^{- \lambda x} dx \right) \mu (d \xi)  \nonumber \\
&=& \frac1{\lambda}\int_{0 }^\infty \left( 1-e^{-\lambda \xi } \right) \mu (d \xi)   
 =  \frac{\phi_0 (\lambda ) }{\lambda}.  \label{e:2.3}
\end{eqnarray}

\medskip

\begin{lemma}\label{L:2.1}
There is a Borel set $\sN \subset (0, \infty)$ having zero Lebesgue measure so that 
  $$
 \P (\bar S_s\geq  t)=\int_0^s \E \left[ w(t-\bar S_r) 1_{\{t\geq \bar S_r \}}\right]  dr
 \quad \hbox{for every } s>0 \hbox{ and } t\in (0, \infty) \setminus \sN .
 $$
 Consequently, for every $t\in (0, \infty) \setminus \sN$, $s\mapsto \P (\bar S_s\geq  t)$ is continuous 
 and $\P (\bar S_s= t)=0$ for every $s>0$. 
 \end{lemma}

\pf  Note that since $r\mapsto \bar S_r$ is strictly increasing a.s.,   by Fubini theorem,
$$
\int_0^s \E \left[ w(t- \bar S_r) 1_{\{t\geq \bar S_r \}}\right]  dr = \int_0^s \E \left[ w(t- \bar S_r) 1_{\{t> \bar S_r \}}\right]  dr.
$$
 For each fixed $s>0$, the Laplace transform of $t\mapsto \P  (\bar S_s \geq t)$ is 
\begin{eqnarray*}
\int_0^\infty e^{-\lambda t}  \P  (\bar S_s \geq t) dt  
&=& \int_0^\infty e^{-\lambda t}  \P  (\bar S_s >t) dt \\
&=&  -
\frac1{\lambda} \int_0^\infty \P (\bar S_s >t) d e^{-\lambda t} 
= \frac1{\lambda} + \frac1{\lambda}  \int_0^\infty e^{-\lambda t} d_t \P (\bar S_s>t) \\
&=& \frac1{\lambda}  - \frac1{\lambda}\E \left[ e^{-\lambda \bar S_s} \right] = \frac{1-e^{-s \phi _0 (\lambda)}}{\lambda}. 
\end{eqnarray*}
By  Fubini theorem and \eqref{e:2.3}, the Laplace transform of $t\mapsto \int_0^s \E \left[ w(t-\bar S_r) 1_{\{t\geq \bar S_r\}}\right]  dr$ is
\begin{eqnarray*}
\int_0^\infty e^{-\lambda t}   \left(\int_0^s \E \left[ w(t-\bar S_r) 1_{\{ t\geq \bar S_r \}}\right]  dr \right) dt 
&=& \int_0^s \E \left[ \int_0^\infty e^{-\lambda t}  w(t-\bar S_r) 1_{\{t > \bar S_r \}} dt \right]  dr  \\
&=& \int_0^s \E \left[ e^{- \lambda \bar S_r} \int_0^\infty e^{-\lambda x  }  w(x)   dx \right]  dr  \\ 
&=& \frac{\phi_0 (\lambda)}
{\lambda} \int_0^s e^{-r \phi_0  (\lambda )} dr 
= \frac{1-e^{-s\phi_0 (\lambda)}}{\lambda}, 
\end{eqnarray*}
which is the same as the Laplace transform of $t\mapsto \P  (\bar S_s >t)$. 
By the uniqueness of the Laplace transform that for each fixed $s>0$, 
\begin{equation} \label{e:2.4}
\P (\bar S_s\geq  t)=\int_0^s \E \left[ w(t-\bar S_r) 1_{\{t\geq \bar S_r \}}\right]  dr 
\end{equation}
for a.e. $t>0$. 
Hence there is a  Borel subset $\sN \subset (0, \infty)$ having zero Lebesgue measure so that
\eqref{e:2.4} holds for every $t\in (0, \infty) \setminus \sN$ and for every rational $s>0$.
Note that   for each fixed $t>0$, $s\mapsto \P (\bar S_s\geq  t)$ is right-continuous.
On the other hand,    
for each fixed $t>0$, 
$s\mapsto \int_0^s \E \left[ w(t- \bar S_r) 1_{\{t\geq \bar S_r \}}\right]  dr $ is  continuous.
It follows that \eqref{e:2.4} holds for every $t\in (0, \infty) \setminus \sN$ and   every   $s>0$. 
 Consequently, 
for every $t\in (0, \infty) \setminus \sN$, $s\mapsto   \P (\bar S_s\geq  t)$ is continuous.
Since the subordinator $t\mapsto \bar S_t$ is strictly increasing a.s. and is stochastically continuous in the sense that
$\P (\bar S_r=\bar S_{r-})=1$ for all $r>0$, we have
$$
\P (\bar S_s\geq  t)= \lim_{r\uparrow s} \P (\bar S_r\geq t) = \P(\bar S_s >t) \quad \hbox{for every } s>0.
$$
In other words, $\P(\bar S_s=t)=0$ for every $t\in (0, \infty) \setminus \sN$ and all $s>0$. 
\qed

\bigskip

Define $G(0)=0$ and $G(x)= \int_{0}^x w(t) dt$ for $x>0$. Then by \eqref{e:2.2}, 
$G(x)$ is a continuous function on $[0, \infty)$
with $G'(x)=w(x)$ on $(0, \infty)$. By the integration by parts formula, for every $t>0$, 
\begin{eqnarray} \label{e:2.5} 
\int_0^t w(t-r) \P (S_s>r) dr
&=& -\int_0^t \P(S_s>r) d_r G(t-r)  \nonumber \\
&=& G(t)+\int_0^t G(t-r) d_r \P(S_s >r) \nonumber \\
&=& G(t)-\int_0^t G(t-r) d_r \P(S_s \leq r) \nonumber \\
&=& G(t)-\E \left[ G(t-S_s) 1_{\{t\geq S_s \}}\right].
\end{eqnarray}
 In particular,  
 $$
\E \left[ G(t-S_s) 1_{\{t\geq S_s \}}\right] \leq G(t)
 \quad \hbox{for every } t>0.
  $$
 For  each fixed $t>0$,   by   \eqref{e:2.2} and dominated convergence theorem,  
 $$
 s\mapsto \int_0^t w(t-r) \P (S_s>r) dr = \int_0^t w(t-r) \P (S_s\geq r) dr
 $$
  is a right continuous increasing function. Hence by \eqref{e:2.5},  
 $s\mapsto \E \left[ G(t-S_s) 1_{\{t\geq S_s \}}\right]$ is a right continuous decreasing   function on $[0, \infty)$.

 \medskip
 
 \begin{cor}\label{C:2.2} Let $\sN\subset (0, \infty)$ be the  set  in Lemma \ref{L:2.1}, which has zero Lebesgue measure.
 \begin{description}
 \item{\rm  (i)} $\displaystyle \int_0^\infty \E \left[ w(t-\bar S_r) 1_{\{t\geq \bar S_r \}}\right]  dr =1$ for every $t\in (0, \infty)\setminus \sN$.
 
 \item{\rm (ii)}  $\displaystyle \int_0^\infty \E \left[ G(t- \bar S_r) 1_{\{ t\geq \bar S_r\}}\right]  dr =t$ for every $t>0$.
 
 \item{\rm (iii)}  $\displaystyle \int_0^\infty \E \left[ G(t- \  S_r) 1_{\{ t\geq   S_r\}}\right]  dr \leq t$ for every $t>0$.
 \end{description}
 \end{cor}
 
 \pf  (i) just follows from Lemma \ref{L:2.1} by taking $s\to \infty$. 
 
 (ii) For $t>0$,  we have by (i) and Fubini theorem that 
 \begin{eqnarray*}
 t &=& \int_0^t  \left( \int_0^\infty \E \left[ w(s- \bar S_r) 1_{\{s\geq 
 \bar S_r 
  \}}\right]  dr \right) ds \\ 
&  =&    \int_0^\infty  \E \left[ \int_0^t  w(s-\bar S_r) 1_{\{ s\geq 
 \bar S_r 
 \}} ds \right]  dr \\
 &=& \int_0^\infty   \E \left[ G(t-\bar S_r) 1_{\{ t\geq \bar S_r \}} \right] dr. 
 \end{eqnarray*}
 
 (iii) Since $G(x)$ is an increasing function in $x$, we have by (ii)
 $$
 \int_0^\infty \E \left[ G(t-    S_r) 1_{\{ t\geq   S_r\}}\right]  dr
 \leq   \int_0^\infty \E \left[ G(t- \bar  S_r) 1_{\{ t\geq   \bar S_r\}}\right]  dr
   =t. 
 $$
This proves the corollary. \qed 

\bigskip

 We define the generalized Caputo derivative $\partial^w_t$ by
\begin{equation}\label{e:caputo}
\partial^w_t f(t):= \frac{d}{dt} \int_0^t w(t-s) (f(s)-f(0)) ds,
\end{equation}
whenever it is well-defined in some function space of $f$. 

\medskip

 Suppose that $\{T_t; t\geq 0\}$ is  a strongly continuous   semigroup with 
infinitesimal generator $(\sL, \sD (\sL))$ in some Banach space $({\mathbb B},  \| \cdot \|)$
 with the property that $ \sup_{t>0} \| T_t  \|  <\infty$.  Here $\| T_t\|$ denotes the operator norm
 of the linear map $T_t: \bB \to \bB$. Note that by the uniform boundedness principle,
  $ \sup_{t>0} \| T_t  \|  <\infty$  is equivalent to $\sup_{t>0} \| T_t f\|  <\infty$ for every $f\in \bB$. 
 Typical examples of such uniformly bounded strongly continuous semigroups are:
 
 \begin{description}
 \item{(i)} Transition semigroup $\{P_t; t\geq 0\}$ of a strong Markov process $X=\{X_t, t\geq 0; \, \P_x, x\in E\}$ on a Lusin  space $E$ 
 that has a weak dual with respect to some reference measure $m$ on $E$. Then for every $p\geq 1$,  $\{P_t; t\geq 0\}$ is a strongly continuous semigroup
 in $\bB:=L^p(E; m)$  with $\sup_{t>0} \| P_t\|_{p\to p} \leq 1$. 
 The infinitesimal generator $(\sL, \sD (\sL))$ of $\{P_t; t\geq 0\}$ in $ L^p(E; m)$ is called the $L^p$ generator of the Markov process $X$.

 \item{(ii)} Transition semigroup $\{P_t; t\geq 0\}$ of a Feller process $X=\{X_t, t\geq 0; \, \P_x, x\in E\}$
   on a locally compact separable Hausdorff space  $E$.
 In this case, $\{P_t; t\geq 0\}$ is a strongly continuous semigroup 
 in the space $(C_\infty (E), \| \cdot \|_\infty)$ of continuous functions  on $E$ that vanish at infinity equipped with uniform norm. 
The infinitesimal generator $(\sL, \sD (\sL))$ of $\{P_t; t\geq 0\}$ in $\bB:=(C_\infty (E), \| \cdot \|_\infty)$ 
 is called the Feller generator of  $X$.  
   
\item{(iii)} Certain Feynman-Kac semigroups (can be non-local Feynman-Kac semigroups or even generalized Feynman-Kac semigroups)
in $L^p$-space or  in $C_\infty (E)$  
of a Hunt process $X$; cf.  \cite{CFKZ, CS}. 
\end{description}

\medskip

For $\alpha >0$, let $G_\alpha:=\int_0^\infty e^{-\alpha t} T_t dt$ be the resolvent of the semigroup $\{T_t; t\geq 0\}$ on Banach space
${\mathbb B}$. Then by the resolvent equation, $\sD (\sL)= G_\alpha ({\mathbb B}) = G_1 ({\mathbb B})$, which is
dense in the Banach space $({\mathbb B}, \| \cdot \|)$. 

\medskip

Let $E_t:= \inf\{s>0: S_s >t\}$, $t\geq 0$, be the inverse subordinator.
Define 
\begin{equation}\label{e:2.7}
 u(t, x)= \E  \left[ T_{E_t} f (x) \right] = \int_0^\infty T_s f(x) d_s \P (E_t \leq s) = \int_0^\infty T_s f(x) d_s \P (S_s \geq t).
\end{equation}

The following is the main result of this paper, which gives the existence and uniqueness of solutions to time fractional equation
\eqref{e:2.8}. Theorem \ref{T:1.1} is its particular case, where $T_t$ is the transition semigroup of a strong Markov process $X$ given by
$T_t f(x)=\E_x [ f(X_t)]$.  

\begin{thm}\label{T:2.3}
Suppose that $(\sL, \sD (\sL))$ is the infinitesimal generator of  a uniformly bounded strongly continuous  semigroup $\{T_t; t\geq 0\}$ 
in a Banach space $({\mathbb B},  \| \cdot \|)$.
 For  every $f\in \sD (\sL)$, $u(t, x):= \E \left[ T_{E_t} f (x) \right]$ is a solution in   $(\bB, \| \cdot \|)$ to
\begin{equation}\label{e:2.8} 
\left( \kappa \partial_t  + \partial^w_t \right) u(t, x) = \sL u(t, x) \quad \hbox{with } u(0, x)=f(x) 
\end{equation}
 in the following sense:

\begin{description}
\item{\rm (i)} $\sup_{t>0} \| u(t, \cdot )\|<\infty$, $x\mapsto u(t, x)  $ is in $\sD (\sL)$ for each $t\geq 0$  
with $\sup_{t\geq  0} \| \sL u(t, \cdot )\| <\infty$, 
and  both $t\mapsto u(t, \cdot)$ and $t\mapsto \sL u(t, \cdot)$ are continuous in $(\bB, \| \cdot \|)$; 

\item{\rm (ii)} for every $t>0$, $I^w_t (u)  :=\int_0^t w(t-s) (u(s, x) -f(x)) ds$ is absolutely convergent in   $(\bB, \| \cdot \|)$ and 
$$
  \lim_{\delta \to 0} \frac{1}{\delta}  \left( \kappa (u(t+\delta, \cdot)-u(t, \cdot) ) + 
   I_{t+\delta}^w(u) - I_t^w(u) \right) = \sL u(t, x) \quad \hbox{in }  (\bB, \| \cdot \|). 
$$  
\end{description}
When $\kappa >0$, $t\mapsto u(t, \cdot)$ is globally  Lipschitz continuous in $ (\bB, \| \cdot \|)$ and hence
$\partial_t u(t, \cdot)$ exists in $ (\bB, \| \cdot \|)$ for a.e. $t\geq 0$. \footnote{See Section \ref{S:4} for an improved statement.}

Conversely, if $u(t, x)$ is a solution to \eqref{e:2.8} in the sense of  {\rm (i)}  and {\rm (ii)}  above with $f\in \sD (\sL)$,
 then $u(t, x)= \E \left[ T_{E_t} f (x) \right]$ in $\bB$
for every $t\geq 0$. 
\end{thm}

\pf   (a) (Existence)
 Clearly  for $f\in \sD (\sL)$, 
$$
\sup_{t>0} \| u(t,  \cdot )\|  \leq \sup_{t>0} \E \left[  \| T_{E_t} f  \| \right] \leq  \sup_{r> 0} \| T_r  f\| <\infty.
$$ 
By the same reason, $\sup_{t>0}   \E \left[  \| T_{E_t}  \sL f  \| \right] \leq \sup_{r> 0} \| T_r  \sL f\| <\infty$. 
 Since
$$
\lim_{\delta \to 0} \frac1{\delta} (T_\delta u(t, \cdot)-u(t. \cdot)) = \lim_{\delta \to 0}  \E \left[ T_{E_t} (T_\delta f -f )/\delta))\right]
= \E \left[ T_{E_t} \sL f \right]
$$
in $(\bB, \| \cdot \|)$. we conclude that  $u(t, \cdot) \in \sD (\sL)$ with $\sL u(t, \cdot) = \E \left[ T_{E_t} \sL f \right]$ for every $t>0$.
 Since $\{T_t; t\geq 0\}$ is a strongly continuous semigroup on $\bB$ with $\sup_{t\geq 0} \| T_t  \| <\infty $ and
$t\mapsto E_t$ is continuous a.s.,   we have by bounded convergence theorem that  both 
$t\mapsto u(t, \cdot) = \E \left[T_{E_t} f  \right]$ 
and $t\mapsto \sL u(t, \cdot) =\E \left[    T_{E_t}  (\sL f ) \right]$ are  continuous in $(\bB, \| \cdot \|)$.

It follows from \eqref{e:2.7}, \eqref{e:2.5},    and the integration by parts formula that for every $t>0$, 
\begin{eqnarray*}
&& \int_0^t w(t-r) (u(r, x)-u(0, x)) dr \\
&=& \int_0^t w(t-r) \left( \int_0^\infty (T_s f(x)-f(x)) d_s \P (S_s \geq r) \right)dr  \\
&=& \int_0^\infty  (T_s f(x)-f(x))  d_s \left( \int_0^t w(t-r) \P (S_s > r) dr \right) \\
&=& - \int_0^\infty  (T_s f(x)-f(x))  d_s\E \left[ G(t-S_s) 1_{\{ t\geq S_s \}}\right]   \\
&=& \int_0^\infty  \E \left[ G(t-S_s) 1_{\{t\geq S_s \}}\right] \sL T_s f(x)  ds .
\end{eqnarray*}
Note that since $\sup_{s>0} \| T_s f\|  <\infty$ and $ \sup_{s>0} \| \sL T_s f\|  =\sup_{s>0} \|  T_s \sL f\|  <\infty$,
by    \eqref{e:2.2}   and Corollary \ref{C:2.2}, all the integrals in above display are absolutely convergent in the Banach space 
$(\bB, \| \cdot \|)$, while the second 
 equality 
 is justified by the Riemann sum approximation of Stieltjes integrals, Fubini theorem
and the dominated convergence theorem.\footnote{Since 
$\| \int_0^\infty  (T_s f-f) d_s \P (S_s>r)\| \leq (M+1) \| f\|$ and  $\int_0^t w(t-r)dr=\int_0^t w(s)ds<\infty$, 
  by  Riemann sum approximation,  the dominated convergence theorem and Fubini's theorem, for partitions $\Pi$ of $[0, \infty)$,
\begin{eqnarray*}
&&\int_0^t w(t-r) \left( \int_0^\infty (T_s f(x)-f(x)) d_s \P (S_s \geq r) \right)dr  \\
&=& \int_0^t w(t-r) \lim_{\| \Pi\|\to 0} \sum_i  (T_{s_i} f(x)  -f(x))(\P (S_{s_{i+1}}>r)-\P (S_{s_i}>r)) dr\\
&=& \lim_{\| \Pi\|\to 0} \int_0^t w(t-r) \sum_i  (T_{s_i} f(x)  -f(x)) (\P (S_{s_{i+1}}>r)-\P (S_{s_i}>r)) dr\\
&=& \lim_{\| \Pi\|\to 0}  \sum_i    (T_{s_i} f(x)  -f(x))  \left( \int_0^t w(t-r) \P (S_{s_{i+1}}>r)dr -\int_0^t w(t-r)  \P (S_{s_i}>r)) dr \right)\\
&=& \int_0^\infty   (T_{s } f(x)  -f(x))   d_s \left( \int_0^t w(t-r) \P (S_s>r) dr \right). 
\end{eqnarray*} 
}
On the other hand,   $\P (S_r \geq s )=1$ when $s\leq \kappa r $, 
while for a.e. $s \in (\kappa r, \infty)$, we have by Lemma \ref{L:2.1} that 
\begin{equation}\label{e:2.9}
\P (S_r \geq s ) = \P (\bar S_r \geq s-\kappa r) =\int_0^r \E \left[ w(s-\kappa r -\bar S_y) 1_{\{s-\kappa r > \bar S_y \}}\right]  dy.
\end{equation}
So for every $t>0$, 
\begin{eqnarray} \label{e:2.10}
\int_0^t \P (S_r \geq s ) ds & =&
 (\kappa r ) \wedge t + \E   \int_0^r    
\left( \int_{ (\kappa r ) \wedge t}^t  w(s-\kappa r -\bar S_y) 1_{\{s-\kappa r > \bar S_y \}}   ds\right) dy  \nonumber \\
 & =&
 (\kappa r ) \wedge t + 1_{\{\kappa r  <t\} } \E   \int_0^r  G(t-\kappa r -\bar S_y) 1_{\{ t-\kappa r > \bar S_y\}} dy . 
\end{eqnarray}
 Since 
\begin{eqnarray*}
\sL u(s, x) &=& \sL \E  \left[ T_{E_s} f (x) \right] = \E  \left[ T_{E_s} \sL f (x) \right] \\
&=& \int_0^\infty   T_r \sL  f (x) d_r \P (E_s \leq r)
= \int_0^\infty   T_r  \sL f (x) d_r \P (S_r \geq s ) ,   
\end{eqnarray*} 
we have  by \eqref{e:2.9} and \eqref{e:2.10} that 
\begin{eqnarray*}
&& \int_0^t \sL u(s, x) ds \\
 &=& \int_0^t \left(\int_0^\infty   T_r  \sL f (x) d_r \P (S_r \geq s ) \right)ds  \\
&=& \int_0^\infty   T_r  \sL f (x) d_r \left( \int_0^t \P (S_r \geq s ) ds \right)   \\
&=& \E \int_0^{t/\kappa} T_r  \sL f (x) \left( \kappa   + G( t-\kappa r -\bar S_r) 1_{\{ t-\kappa r > \bar S_r\}} -\kappa \int_0^r 
w(  t-\kappa r -\bar S_y)    1_{\{ t-\kappa r > \bar S_y\}}   dy \right)   dr   \\
 &=& \int_0^\infty  T_r \sL f(x)  \E \left[ G(t-S_r) 1_{\{t\geq S_r \}} \right] dr
 + \kappa   \int_0^{t/\kappa} T_r  \sL f (x)  \left(1-\P(S_r \geq t) \right)   dr \\
  &=& \int_0^\infty  T_r \sL f(x)  \E \left[ G(t-S_r) 1_{\{t\geq  S_r \}} \right] dr
 + \kappa \int_0^\infty  \P(S_r  <  t) d_r \left( T_r    f (x)   -f(x)\right) \\
  &=& \int_0^\infty  T_r \sL f(x)  \E \left[ G(t-S_r) 1_{\{t\geq S_r \}} \right] dr
 + \kappa \int_0^{t/\kappa} \P (E_t > r)  d_r \left( T_r    f (x)   -f(x)\right) \\
  &=& \int_0^\infty  T_r \sL f(x)  \E \left[ G(t-S_r) 1_{\{t\geq S_r \}} \right] dr
 + \kappa \int_0^\infty  \left( T_r    f (x)   -f(x)\right)  d_r \P (E_t\leq r)   \\
   &=& \int_0^\infty  T_r \sL f(x)  \E \left[ G(t-S_r) 1_{\{t\geq S_r \}} \right] dr
 + \kappa \E \left[ T_{E_t} f(x)    -f (x) \right] \\
 &=& \int_0^\infty  T_r \sL f(x)  \E \left[ G(t-S_r) 1_{\{t\geq S_r \}} \right] dr
 + \kappa (u(t, x)-u(0, x)). 
 \end{eqnarray*}
  Thus we have for every $t>0$, 
$$ 
\kappa (u(t, x)-u(0, x)) + \int_0^t w(t-r) (u(r, x)-u(0, x)) dr = \int_0^t \sL u(s, x) ds.
$$
Consequently,  $\left( \kappa \partial_t +\partial^w_t \right) u(t, x) = \sL u(t, x)$ in $\bB$ as $t\mapsto \sL u(t, \cdot)$ is continuous in $(\bB, \| \cdot \|)$.  

Since $\{T_t; t\geq 0\}$ is a uniformly bounded strongly continuous semigroup in  $(\bB, \| \cdot \|)$, for $f\in \sD (\sL)$ and $t_2>t_1\geq 0$, 
$$
\| T_{t_2} f - T_{t_1} f \| \leq  \int_{t_1}^{t_2} \| \partial_s T_s f \| ds = \int_{t_1}^{t_2} \| \sL T_s f \| ds 
= \int_{t_1}^{t_2} \| \  T_s \sL f \| ds \leq c \| \sL f \| \, |t_2-t_1|  .
$$
Note that when $\kappa >0$,  $|E_t-E_s| \leq |t-s|/\kappa$. Hence we have from the above display that for every $t>s\geq 0$, 
$$
\| u (t, \cdot ) - u(s, \cdot)\| = \| \E \left[ T_{E_t} f  -T_{E_s} f \right]  \| 
\leq c_1 \E |E_t-E_s | \leq c_2 (t-s);
$$
that is, $t\mapsto u(t, \cdot)$ is globally Lipschitz continuous in $(\bB, \| \cdot \|)$.
This implies in particular  that $u(t, \cdot)$ is differentiable in $t$ as an element in $(\bB, \| \cdot \|)$
for a.e. $t>0$. 

\medskip

(b) (Uniqueness) Suppose that $u(t, x)$ is a solution to \eqref{e:2.8} in the sense of (i) and (ii) with $f\in \sD (\sL)$.
Then $v(t, x):= u(t, x)- \E \left[ T_{E_t} f (x) \right]$ is a solution to \eqref{e:2.8} with $v(0, x)=0$. 
Hence we have for every $t>0$, 
\begin{equation}\label{e:2.11}
\kappa v(t, x)+\int_0^t w(t-r) v(r, x) dr = \int_0^t \sL v(s, x) ds.
\end{equation}
Let $V(\lambda, x):=\int_0^\infty e^{-\lambda t} v(t, x) dt$, $\lambda >0$,
 be the Laplace transform of $t\mapsto v(t, x)$. Clearly for every $\lambda >0$, $V(\lambda, \cdot)\in \bB$
with $\| V(\lambda, \cdot)\| \leq \lambda^{-1} \sup_{t>0} \| v(t, \cdot )\|$. 
 Since $v(t, \cdot) \in \sD (\sL)$ for every $t>0$ with $\sup_{t>0}\  \| \sL v(t, \cdot) \|<\infty$, 
we have by dominated convergence theorem that for every $\lambda >0$, 
\begin{eqnarray*}
\lim_{\delta \to 0} \frac1{\delta} (T_\delta V(\lambda, \cdot) -V(\lambda, \cdot)) 
&=& \lim_{\delta \to 0} \int_0^\infty e^{-\lambda t} \frac1{\delta} (T_\delta v(t, \cdot) -v(t, \cdot)) dt  \\
&=& \lim_{\delta \to 0} \int_0^\infty e^{-\lambda t}  \left(\frac1{\delta} \int_0^\delta T_s \sL v(t, \cdot) ds \right) dt \\
&=&  \int_0^\infty e^{-\lambda t} \sL v(t, \cdot) dt. 
\end{eqnarray*}
This shows that for 
each $\lambda >0$, $V(\lambda, \cdot )\in \sD (\sL)$ with  
$$ 
\sL V(\lambda, \cdot) = \int_0^\infty e^{-\lambda t} \sL v(t, \cdot) dt
\quad \hbox{and} \quad 
\|\sL V(\lambda, \cdot) \|\leq \int_0^\infty e^{-\lambda t}  \| \sL v(t, \cdot) \| dt\leq \frac1{\lambda} \sup_{t>0} \| \sL v(t, \cdot)\|.
$$
Taking Laplace transform in $t$ on both sides of \eqref{e:2.11} yields 
$$
V(\lambda, x) \left( \kappa + 
 \int_0^\infty e^{-\lambda s} w(s) ds \right)
= \frac{1}{\lambda} \int_0^\infty  e^{-\lambda t} \sL v(t, x) dt = \frac{\sL V(\lambda, x)}{\lambda}.
$$ 
Thus by \eqref{e:2.1} and \eqref{e:2.3}, $\sL V(\lambda, x) = \left(\kappa \lambda + \phi_0 (\lambda) \right) V(\lambda, x)
=\phi (\lambda ) V(\lambda, x)$. In other words, 
$$ 
(\phi (\lambda)-\sL ) V(\lambda, x) = 0\quad \hbox{for every } \lambda >0.
$$
Since $\sL$ is the infinitesimal generator of a uniformly bounded strongly continuous semigroup $\{T_t, t\geq 0\}$ in Banach space $\bB$,
for every $\alpha >0$, the resolvent $G_\alpha = \int_0^\infty e^{-\alpha t} T_t dt$ is well defined and is the inverse to $\alpha-\sL$.
Hence we have from the last display that $V(\lambda, \cdot)=0$ in $\bB$ for every $\lambda >0$. By the uniqueness of Laplace transform,
we have $v(t, \cdot)=0$ in $\bB$ for every $t>0$. This establishes that $u(t, x)= \E \left[ T_{E_t} f  (x)\right]$ in $\bB$ for every $t\geq 0$. 
\qed

\begin{remark}\label{R:2.4} \rm \begin{description}
\item{(i)} The assumption that $f  \in \sD (\sL)$   in Theorem \ref{T:2.3}
is to ensure that all the integrals involved in the proof of Theorem \ref{T:2.3} are absolutely convergent in the Banach space $\bB$.
This condition can be relaxed  if we formulate the equation \eqref{e:2.8} in the weak sense when the uniformly bounded 
strongly continuous semigroup
$\{T_t;t\geq 0\}$ is symmetric in  a Hilbert space $L^2(E; m)$ and so its quadratic form can be used to formulate weak solutions. 
This will be carried out in  the ongoing joint work  \cite{CKKW} with  Kim,   Kumagai and   Wang.  
It in particular applies to the case where $\{T_t;  t\geq 0\}$ is the transition semigroup of any $m$-symmetric
Markov process on a Lusin  space $E$, which is a strongly continuous contraction symmetric semigroup in $L^2(E; m)$.

\item{(ii)} There are two closely related  work    \cite{MS2, Ko}.
Suppose that $X=\{X_t,  t\geq 0; \P_x, x\in \bR^d\}$ is a L\'evy process on $\bR^d$   and generator $\sL$,  and 
$S=\{S_t; t\geq 0\}$ is a driftless subordinator 
with Laplace exponent $\phi$ and L\'evy measure $\mu$. 
Let $E_t:=\inf\{s>0: S_s >t\}$ be the inverse subordinator. 
Under the assumption that $\kappa =0$, $\mu (0, \infty) =\infty$,
 $\int_0^1 x | \log x| \mu (dx)<\infty$ and that the L\'evy process $X$ has a transition density function, it is shown in
  \cite[Theorem 4.1]{MS2} that $u(t, x):= \E_x \left[ f( X_{E_t})\right]$ is a mild solution of 
  the following pseudo-differential equation
  $$ 
  \phi (\partial_t) u (t, x) = \sL u (t, x) + f(x) \mu (t, \infty) . 
  $$
Here $\phi (\partial_t)$ is a pseudo-differential operator in time variable $t$ formulated using Fourier multiplier.

Under the assumption that the L\'evy measure $\mu  $ of the subordinator $S_t$ satisfying condition
$\mu (d\xi )\geq \xi^{1+\beta} d\xi$ on $(0, \eps)$ for some $\eps >0$ and $\beta >0$, and $\{T_t; t\geq 0\}$
is the transition semigroup of a Feller process $X=\{X_t, t\geq 0; \P_x, x\in \bR^d\}$
 on $\bR^d$  whose domain of infinitesimal generator 
contains $C^2(\bR^d) \cap C_\infty (\bR^d)$, 
\cite[Theorem 8.4.2]{Ko} asserts that for every $f\in C^2(\bR^d) \cap C_\infty (\bR^d)$, 
$u(t, x):=\E_x \left[ f(X_{E_t})\right]$ satisfies 
$$
A_t^* u(t, x)= \sL u(t, x) + f(x) A^* (1_{(0, \infty)}) (t) \quad \hbox{with } u(0, x)=f(x), 
$$
where $A^*$ is the dual of the infinitesimal generator of the subordinator $S_t$ and notation $A^*_t u(t, x)$ means
that the operator $A^*$ is applied to the function $t\mapsto u(t, x)$.
Here $C^2(\bR^d)$ is the space of $C^2$-smooth functions on $\bR^d$ and $C_\infty (\bR^d)$ is the space
of continuous functions on $\bR^d$ that vanish at infinity. 
In \cite[Theorem 8.4.2]{Ko} , the subordinator $S_t$ may have drift $\kappa \geq 0$. 

Similar problem has also been considered in  \cite{T} under more restrictive conditions and using a different approach.
The time fractional derivative there is of the form 
$$
\kappa  \frac{\partial  u(t)}{\partial t}+ \int_0^t  w(s) \frac{\partial}{\partial t} u(t-s)  ds   .
$$
This requires  regularity assumption beyond absolute continuity on the function $t\mapsto u(t)$, as $w(s)$ is unbounded near $s=0$.
The absolute convergence
of the singular integral should be checked and   justified. 

\item{(iii)} Suppose the subordinator $S$ is driftless and has L\'evy measure 
$\mu (dx)= \left( \int_0^1 \frac{\beta}{x^{1+\beta}} \frac{c(\beta)}{\Gamma (1-\beta)} d\beta \right) dx$,
where $c(\beta)\geq 0$ is a measurable function with $\int_0^1 c(\beta) d\beta <\infty$.
(Note that $\Gamma (1-\beta) \asymp \frac{1}{1-\beta}$ for $0<\beta <1$.)   
Then $w(x):=\mu (x, \infty) = \int_0^1 x^{-\beta}\, \frac{c(\beta)}{\Gamma (1-\beta)} d\beta$.
The time fractional derivative $\partial^w_t$ defined in this paper is the 
distributed-order fractional derivative defined in \cite{MS1}. In this case, for continuously differentiable function $f$ on $[0, \infty)$, the time fractional derivative $\partial^w_t f(t)$
is the mixture of Caputo derivatives of order $\beta$'s:   
$$ \partial^w_t f (t) = \int_0^1 \partial^\beta_t f(t) \, c(\beta) d\beta.
$$

\item{(iv)} Cauchy problems with distributed order time fractional derivatives (where $\kappa=0$) 
were also studied in \cite{MNV2} 
for uniformly elliptic generators of divergence form
in bounded $C^{1, \gamma}$  domains with Dirichlet boundary condition, under certain regularity conditions of the diffusion matrices. 
We also mention \cite[Theorem 2]{MSc}  where $\{S_t; \geq 0\}$ is a subordinator without drift 
and $\{T_t; t\geq 0\}$ is the transition semigroup of  a one-dimensional diffusion  killed
at certain rate via Feynman-Kac transform.  

\item{(v)} There are  limited results in literature on  the uniqueness for the time fractional equations \eqref{e:2.8}; 
 see \cite{Ko1, Ko2,  MNV1} for cases of $\partial^\beta_t u =\sL u$ and \cite{Lu} for distributed order time fractional
equation $\partial^w_t u =\sL u$ where $\sL$ is a one-dimensional differential operator in a bounded interval.  
   We mention that Remark 3.1 of a recent preprint \cite{BLM} contains a uniqueness result for solutions to 
$\partial^\beta_t u =\sL u$, where
$\sL$ is the Feller generator of a doubly Feller process killed upon leaving a bounded regular domain, proved also by using 
Laplace transform similar to our   uniqueness proof for Theorem \ref{T:2.3} in this paper.  

\item{(vi)} When the uniformly bounded strongly continuous semigroup $\{T_t; t\geq 0\} $ in Theorem \ref{T:2.3} 
has an integral kernel $p(t, x, y)$ with  respect to some measure $m(dx)$, then there is a kernel $q(t, x, y)$ so that
$$
u(t, x):= \E  \left[ T_{E_t} f(x) \right] =\int_E q(t, x, y) f(y) m(dy);
$$
in other words, 
$$
q(t, x, y):= \E \left[ p(E_t, x, y) \right] = \int_0^\infty p(s,x, y) d_s \P (E_t \leq s)
$$
 is the fundamental solution to the time fractional equation
 $\left( \kappa \partial_t+ \partial^w_t\right) u = \sL u$ under the setting of this paper. In \cite{CKKW}, two-sided estimates on $q(t, x, y)$ are obtained
 when $\kappa =0$ and $\{T_t; t\geq 0\}$ is the transition semigroup 
of a diffusion process that satisfies two-sided Gaussian-type estimates or of a stable-like process on metric measure spaces.  
\end{description}
\end{remark}

\begin{example} \label{E:2.5}
\rm \begin{description} \item{(i)} When $\{S_t; t\geq 0\}$ is a  
 $\beta$-stable subordinator with $0<\beta<1$ with
Laplace exponent $\phi (\lambda) = \lambda^\beta$, 
it is easy to check that  $S_t$ has no drift (i.e. $\kappa=0$)  and its L\'evy measure is $\mu (dx)= \frac{\beta}{\Gamma (1-\beta)} x^{-(1+\beta)} dx$. 
Hence
$$
w(x):= \mu (x, \infty)= \int_x^\infty \frac{\beta}{\Gamma (1-\beta)} y^{-(1+\beta)} dy = \frac{x^{-\beta}}{\Gamma (1-\beta)}.
$$
Thus  the time fractional derivative $\partial^w_t f $ defined by \eqref{e:1.2} is exactly the Caputo derivative of order $\beta$ defined by \eqref{e:1.1}. In this case, Theorem \ref{T:2.3} recovers the main result of \cite{BMe} and \cite[Theorem 5.1]{MS0}. 

\item{(ii)} We call a subordinator $\{S_t; t\geq 0\}$  truncated $\beta$-stable 
subordinator if  it is driftless and its L\'evy measure   is 
$$
 \mu_\delta (dx) = \frac{\beta}{\Gamma (1-\beta)}   \,  x^{-(1+\beta)} 1_{(0, \delta] } (x) dx
$$
for some $\delta >0$.  In this case,
$$  w_\delta(x):=\mu_\delta (x, \infty) = 1_{\{0<x\leq \delta\}} \int_x^\delta \frac{\beta}{\Gamma (1-\beta)} y^{-(1+\beta)} dy 
= \frac{1}{\Gamma (1-\beta)} \left( x^{-\beta} - \delta^{-\beta}\right) 1_{(0, \delta]} (x). 
$$
So the corresponding  the fractional derivative   of \eqref{e:1.2} is 
$$
\partial^{w_\delta} _t f (t):= \frac{1}{\Gamma (1-\beta)}  \frac{d}{dt} 
\int_{(t-\delta)^+}^t   \left( (t-s)^{- \beta } -\delta^{-\beta} \right) (f(s)-f(0)) ds . 
$$
This is the fractional-time derivative whose value at time $t$  depends only on  the $\delta$-range of the past of $f$
as mentioned in the Introduction. Theorem \ref{T:2.3} says that the corresponding time fractional equation \eqref{e:1.5}
can be solved by using the inverse of truncated $\beta$-stable subordinator. 
Clearly, as $ \lim_{\delta\to \infty} w_\delta (x) = w(x):=\frac{1}{\Gamma (1-\beta)} x^{-\beta}$. Consequently, the fractional derivative
$\partial^{w_\delta}_t f(t) \to \partial^w_t f(t))$, the  Caputo derivative of $f$ of order $\beta$, in the distributional sense
as $\delta \to 0$.  
Using the probabilistic representation in Theorem \ref{T:2.3}, one can deduce that as $\delta \to \infty$, the solution
to the equation $\partial_t^{w_\delta} u= \sL u $  with  $ u(0, x)=f(x) $
converges to the solution of   $\partial_t^{\beta} u= \sL u $  with  $ u(0, x)=f(x) $. 

If we define 
$$\eta_\delta (r)= \frac{\Gamma (2-\beta) \, \delta^{\beta-1}  }{\beta}w_\delta (r)
=    (1-\beta) \delta^{\beta-1}  \left( x^{-\beta} - \delta^{-\beta}\right) 1_{(0, \delta]} (x), 
$$
then $\eta_\delta (r)$ converges weakly to the Dirac measure concentrated at $0$ as $\delta \to 0$.
So the fractional derivative $\partial_t^{\eta_\delta} f (t)$ converges to $f'(t)$ for every differentiable $f$. 
It can be shown that the subordinator corresponding to $\eta_\delta$, that is,  subordinator with L\'evy measure
$$ 
\nu_\delta (dx) :=\frac{  (1-\beta) \, \delta^{\beta-1}}{\beta  }   \,  x^{-(1+\beta)} 1_{(0, \delta] } (x) dx, 
$$ 
converges as $\delta \to 0$ to deterministic motion $t$ moving at constant speed 1.
Using Theorem \ref{T:2.3}, one can show that the solution
to the equation $\partial_t^{\eta_\delta} u (t, x) = \sL u (t, x)$  with  $ u(0, x)=f(x) $
converges to the solution of the heat equation  $\partial_t  u= \sL u $  with  $ u(0, x)=f(x) $. 
    \qed
\end{description}
\end{example}

\section{Occupation measure for processes time-changed by inverse subordinator}

Suppose $X=\{X_t, t\geq 0; \, \P_x, x\in E\}$ is a general strong Markov process on state space $E$ with infinitesimal generator 
$\sL$, and
$S=\{S_t; t\geq 0\}$ is a subordinator independent of $X$
whose L\'evy measure $\mu$ satisfies $\mu (0, \infty)=\infty$.
Let $\phi$ be the Laplace exponent of $S$; that is, $\E e^{-\lambda S_t} =e^{-t \phi (\lambda)}$.
Note that $\E \left[ S_t \right]= t \phi'(0)$ so in particular $\phi '(0) = \E \left[ S_1 \right]$.  
Let $E_t:=\inf\{s>0: S_s >t\}$ be the inverse subordinator, and $X^*_t:=X_{E_t}$.
Suppose $D$ is an open subset of $E$ and define $\tau_D:=\{t>0: X_t\notin D\}$ to be the first exit time
from $D$ by the process $X$. 
In general, the time-changed process $X^*$ is not a Markov process but we can still define its 
first exit time from $D$ by 
$$
\tau^*_D   := \inf \left\{t>0: X^*_t\notin D \right\}.
$$
Let $\partial$ be a cemetery point. 
The process $X^{*, D}$ defined by $X^{*, D}_t:= X^*_t$ when $t<\tau^*_D$ and $X^*_t:= \partial $
for $t\geq \tau^*_D$ is called the part process of $X^*$ in $D$. 
The part process $X^D$ of $X$ in $D$ is defined in an analogous way. 
We use $\E_x$ to denote mathematical
expectation  taken with respect to
the probability law $\P_x$, under which the Markov process $X$ starts from $x\in E$. 
For every $x\in D$, the occupation measures for $X^D$ and $X^{*, D}$ are defined by
$$
\nu^D_x (A) =\E_x \left[ \int_0^{\tau_D} 1_A (X_s) ds \right] \quad
\hbox{and} \quad  \nu^{*, D}_x (A) =\E_x \left[ \int_0^{\tau^*_D} 1_A (X^*_s) ds \right],
\qquad A\subset D.
$$
Occupation measures describe the average amount of time spent by the processes in subsets of the state space. 

The next theorem says that the occupation measure for the part process $X^{*, D}$ of $X^*$ in $D$
is proportional to that of the part process $X^D$ of $X$ in $D$ when $\phi '(0) <\infty$, 
that is, when the subordinator $S_t$ has finite
mean. When the subordination $S_t$ has infinite mean,  the occupation measure for the
part process $X^{*, D}$ of $X^*$ in $D$ is always infinite.

\begin{thm}\label{T:3.1} For every measurable function $f\geq 0$ on $D$ and $x\in D$,
$$
\E_x  \left[\int_0^{\tau^*_D} f(X^*_t ) dt \right] = \, \phi'(0)  \, \E_x \left[ \int_0^{\tau_D} f(X_t)d t  \right]
  = \, \phi'(0) G_D f(x).
$$
In other words, $\nu_x^{*, D} = \phi'(0) \nu^D_x$ for every open set $D\subset E$ and   every $x\in D$. 
\end{thm}

\pf First note that 
\begin{eqnarray*}
\tau^*_D & =&\inf\{t>0: X_{E_t} \notin D\} 
= \inf\{t>0: E_t=\tau_D\}  \\
&=& \inf\{t>0: S_{\tau_D} >t\} = S_{\tau_D}.
\end{eqnarray*}
For any $f\geq 0$ on $D$, we have using the independence between the strong Markov process $X$ and the subordinator $S$ that 
\begin{eqnarray*}
\E_x  \left[ \int_0^{\tau^*_D} f(X^*_t ) dt \right] 
&=& \E_x \left[ \int_0^{S_{\tau_D}} f(X_{E_t} )dt  \right] = \E_x \left[ \int_0^{\tau_D} f(X_r)d S_r \right]\\
&=& \E_x  \left[ \int_0^{S_{\tau_D}} f(X_{E_t} )dt \right] = \E_x  \E_x \left[   \int_0^{\tau_D} f(X_r)d S_r  \Big| X \right] \\
&=& \E_x   \left[   \int_0^{\tau_D} f(X_r)d (\E S_r)  \right] = \E_x   \left[   \int_0^{\tau_D} f(X_r)d r  \right] \phi'(0)\\
&=& \phi'(0) G_D f(x). 
\end{eqnarray*}
\qed

\begin{remark} \label{R:3.2} \rm   
\begin{description}
\item{(i)}
Taking $f=1$ in Theorem \ref{T:3.1} in particular yields the following relation on mean exit times:  
\begin{equation} \label{e:3.1}
\E_x \left[ \tau^*_D \right] = \phi'(0) \, \E_x \left[ \tau_D \right] \quad \hbox{for every } x\in D.
\end{equation}
When $X$ is either a diffusion process determined by a stochastic differential equation driven by Brownian motion or a rotationally symmetric
$\alpha$-stable process on $\bR^d$, and $\{S_t; t\geq 0\}$ is a tempered  $\beta$-stable subordinator
having Laplace exponent $\phi (\lambda)= (\lambda + m)^\beta -m^\beta$ for some $m>0$ and $0<\beta <1$,  
\eqref{e:3.1} recovers the main result of \cite{DWW}, derived there using  a PDE method.

\item{(ii)} Observe that the part process $X^D$ of $X$ killed upon leaving $D$ is a strong Markov process  in $D$ whose 
infinitesimal generator  $\sL^D$ is $\sL$  in $D$ having zero exterior condition.
The transition semigroup of $X^D$ is  $P^D_tf:=\E_x \left[ f(X_t); t<\tau_D\right]$. 
Hence by Theorem \ref{T:2.3},  for $f\in \sD (\sL^D)$, 
$$
u(t, x):=\E_x \left[ f(X^D_{E_t})\right] =\E_x \left[ f(X^*_t); t<\tau^*_D\right]
$$
is the strong solution to  
$$
\left( \kappa \partial_t  + \partial^w_t \right) u(t, x) = \sL^D u(t, x) \quad \hbox{with } u(0, x)=f(x) \hbox{ in } D.
$$
On the other hand, $ G_D f(x)$ is the solution to the Poisson equation $\sL v=-f$ in $D$ with  $v=0$ on $D^c$.
Hence it follows from Theorem \ref{T:3.1} that for $f\in \sD (\sL^D)$,
$$ G^*_D f(x)=\int_0^\infty \E_x \left[ f(X_t); t<\tau^*_D\right] dt =\int_0^t u(t, x) dt
$$
is the solution to the Poisson equation
$$
\sL^D v=-\phi'(0) f \quad \hbox{in } D \quad \hbox{with } v=0 \hbox{ on } D^c. 
$$ 
Since by (ii) of Theorem \ref{T:2.3} that
$$
  \lim_{\delta \to 0} \frac{1}{\delta}  \left( \kappa (u(t+\delta, \cdot)-u(t, \cdot) ) + 
   I_{t+\delta}^w(u) - I_t^w(u) \right) = \sL u(t, x) \quad \hbox{in }  (\bB, \| \cdot \|). 
$$  
and $t\to \mapsto \sL u(t, \cdot)$ is continuous in $ (\bB, \| \cdot \|)$, we conclude that 
 $$
 \partial_t u(t, \cdot):=\lim_{\delta \to 0} \frac{1}{\delta}  \left(   (u(t+\delta, \cdot)-u(t, \cdot) ) \right)
 = \sL u(t, \cdot) -\partial^w_t u
 $$
 exists and $t\mapsto \partial_t u$ is continuous in $(\bB, \| \cdot \|)$. 
\qed

\end{description}
\end{remark}

\medskip

{\bf Acknowledgement.} The author thanks M. M. Meerschaert for the invitation to the Workshop 
``Future Directions in Fractional Calculus Research and Applications" held  at Michigan State University,
East Lansing,  from October 17-21,  2016, and for helpful comments.
He also thanks T. Kumagai  and J. Wang for helpful comments. 

\section{Note added after publication}\label{S:4}

The last sentence in the existence part of Theorem \ref{T:2.3} can be strengthened as follows:

\medskip

{\it  When $\kappa >0$, $t\mapsto u(t, \cdot)$ is globally  Lipschitz continuous in $ (\bB, \| \cdot \|)$, and 
both $\partial_t u(t, \cdot)$ and $\partial^w_t u (t, \cdot):=\frac{d}{dt} I^w_t(u)$ exists as a continuous function taking values in $ (\bB, \| \cdot \|)$.} 

\medskip

\noindent{\bf Proof.} 
Recall that $\mu$ is the L\'evy measure for the subordinator $S_t$ and $w(x):=\mu (x, \infty)$ for $x>0$.
By \eqref{e:2.2} and the monotone convergence theorem, 
\begin{equation}\label{e:4.1} 
\lim_{r\to 0} rw(r)\leq \lim_{r\to 0} \int_0^r  w(s)ds  =0 . 
\end{equation}  
When $\kappa >0$, it is shown in Theorem \ref{T:2.3} that 
 $t\mapsto u(t, \cdot)$ is globally  Lipschitz continuous in $ (\bB, \| \cdot \|)$.
So the following function taking values in $\bB$ is well defined: 
\begin{equation}\label{e:4.2} 
g(s, \cdot):= (u(s, \cdot)-u(0, \cdot)) w(s) + \int_{(0, s]} (u(s, \cdot)-u(s-r, \cdot)) \mu (dr), \quad s>0,
\end{equation}
with $\lim_{s\to 0} g(s, \cdot )=0$ in $(\bB, \| \cdot \|)$.\footnote{This  definition is motivated by the fact that were $u(r, \cdot)$ differentiable in $r$,  then 
$$ g(s, \cdot) = \int_0^s w(s-r) d_r (u(r, \cdot) -u(s, \cdot)) = \int_0^s w(s-r) d_r  u(r, \cdot) 
$$ 
by  an integration by parts.}    Moreover,  denoting the global Lipschitz constant of $t\to u(t, \cdot)$ in $(\bB, \| \cdot \|)$
by $M_1$, 
we have  for any $s>0$ and $\delta>0$,  
 \begin{eqnarray*}
\| g(s+\delta , \cdot) -g(s, \cdot) \| 
&\leq &   \| u(s+\delta   , \cdot) - u(s,  \cdot) \|  \, w(s+\delta) + \| u(s ,  , \cdot) - u(0,  \cdot) \| \, | w(s+\delta) -w(s)|\\
& & +  \int_{(0, s]}  \|  (u(s+\delta   , \cdot) - u(s+\delta -r   , \cdot) ) -  ( u(s,  \cdot) - u(s-r,  \cdot) \| \, \mu (dr)  \\
&& + \int_{(s, s+\delta]} \|    u(s+\delta   , \cdot) - u(s+\delta -r,  \cdot) \| \,  \mu (dr)   \\
&\leq &  M_1 \delta w(s+\delta) + M_1 s \mu (s, s+\delta]   + M_1   \int_{(s, s+\delta]} r \mu (dr)   \\
&&  +  \int_{(0, s]} \|  (u(s+\delta   , \cdot) - u(s+\delta -r   , \cdot) ) -  ( u(s,  \cdot) - u(s-r,  \cdot) \| \, \mu (dr)     .
\end{eqnarray*} 
In view of \eqref{e:2.2} and \eqref{e:4.1}, each of the first three terms converges to 0 as $\delta \to 0+$, and so does
the fourth term by the dominated convergence theorem and the bound 
$$
\|  (u(s+\delta   , \cdot) - u(s+\delta -r   , \cdot) ) -  ( u(s,  \cdot) - u(s-r,  \cdot) \| \leq 2M_1 r.
$$
Note that  by definition \eqref{e:4.2}, 
$$
g(s, \cdot):= (u(s, \cdot)-u(0, \cdot)) w(s-) + \int_{(0, s)} (u(s, \cdot)-u(s-r, \cdot)) \mu (dr) \quad \hbox{for } s>0.
$$
Using this expression, by a similar argument as above we have for every $s>0$, 
\begin{eqnarray*}
&& \lim_{\delta \to 0+} \| g(s-\delta , \cdot) -g(s, \cdot) \| \\
&\leq & \lim_{\delta \to 0+} \left(  M_1 \delta w(s ) + M_1 s \mu (s-\delta , s)   + M_1   \int_{(s-\delta, s)} r \mu (dr) \right)  \\
&&  +  \lim_{\delta \to 0+} \int_0^{s } 1_{(0, s-\delta]}(r) \|  (u(s   , \cdot) - u(s -r   , \cdot) ) -  ( u(s-\delta ,  \cdot) - u(s-\delta -r,  \cdot) \| \, \mu (dr)   \\
&=& 0  .
\end{eqnarray*}
This establishes the claim that $s\mapsto g(s, \cdot)$ is   continuous in $ (\bB, \| \cdot \|)$ on $[0, \infty)$ with $g(0, \cdot):=0$ .
 For every $t>0$,   
\begin{eqnarray*}
\int_0^t g(s, \cdot) ds 
&=& \int_0^t  (u(s, \cdot)-u(0, \cdot)) w(t-s) ds  +  \int_0^t  (u(s, \cdot)-u(0, \cdot)) (w(s)-w(t-s)) ds \\
&& + \int_0^t \left( \int_{(0, s]} (u(s, \cdot)-u(s-r, \cdot)) \mu (dr) \right)  ds \\
&=& \int_0^t  (u(s, \cdot)-u(0, \cdot)) w(t-s) ds  +  \int_0^t  (u(s, \cdot)-u(t-s,  \cdot)) w(s) ds  \\
&& + \int_{(0, t]}  \left( \int_r^t (u(s, \cdot)-u(s-r, \cdot)) ds \right) \mu (dr)  
\qquad \hbox{ (by Fubini's theorem)} \\
&=&\int_0^t  (u(s, \cdot)-u(0, \cdot)) w(t-s) ds  +  \int_0^t  (u(s, \cdot)-u(t-s,  \cdot)) w(s) ds   \\
&& + \int_{(0, t]} \left( \int_{t-r}^t u(s, \cdot) ds-  \int_0^r u(s, \cdot ) ds \right) \mu (dr) \\
&=&\int_0^t  (u(s, \cdot)-u(0, \cdot)) w(t-s) ds  +  \int_0^t  (u(s, \cdot)-u(t-s,  \cdot)) w(s) ds  ds \\
&& + \int_{(0, t]} \left( \int_0^r  \left( u(t-s , \cdot) - u(s, \cdot ) \right) ds \right) \mu (dr) \\
&=&\int_0^t  (u(s, \cdot)-u(0, \cdot)) w(t-s) ds  +  \int_0^t  (u(s, \cdot)-u(t-s,  \cdot)) w(s) ds  \\
&& + \int_0^t \left( \int_{(s, t]}    \left( u(t-s , \cdot) - u(s, \cdot ) \right) \mu (dr) \right)  ds
\qquad \hbox{ (by Fubini's theorem)} \\
&=&\int_0^t  (u(s, \cdot)-u(0, \cdot)) w(t-s) ds  +  \int_0^t  (u(s, \cdot)-u(t-s,  \cdot)) w(s) ds  \\
&& + \int_0^t   ( u(t-s , \cdot) - u(s, \cdot )  ) (w(s)-w(t))   ds \\
&=&\int_0^t  (u(s, \cdot)-u(0, \cdot)) w(t-s) ds .
\end{eqnarray*}
 Since $t\mapsto g(t, \cdot)$ is continuous in $(\bB, \| \cdot \|)$ over $[0, \infty)$, 
$\partial_t^w u (t, \cdot):=\frac{d}{dt} I^w_t(u) $ exists for every $t>0$
and $t\mapsto \partial_t^w u = g(t, \cdot)$ is continuous in $(\bB, \| \cdot \|)$.
Now it follows from (ii) of Theorem \ref{T:2.3} that
$$
\partial_t u(t, \cdot ):=
  \lim_{\delta \to 0} \frac{1}{\delta}  \left(   (u(t+\delta, \cdot)-u(t, \cdot) )  \right)
  =\frac{1}{\kappa} \left(\sL u(t, \cdot) - \partial^w_t u(t, \cdot) \right) 
  =\frac{1}{\kappa} \left(\sL u(t, \cdot) - g(t, \cdot) \right)
 $$  
exists and $t\mapsto \partial_t u(t, \cdot)$ is continuous in $(\bB, \| \cdot \|)$.
Hence $u(t, x):=\bE \left[ T_{E_t} f(x)\right]$ satisfies
$$ \kappa \partial_t u (t , x) ) + \partial^w_t u (t , x)  = \sL u(t,x) 
$$
in the strong sense in the Banach space $(\bB, \| \cdot \|)$. 
\qed

\bigskip

{\bf Zhen-Qing Chen}

Department of Mathematics, University of Washington, Seattle,
WA 98195, USA

E-mail: \texttt{zqchen@uw.edu}

\end{document}